
%
%
%

\input amstex
\documentstyle{amsppt}
\tolerance=5000

\mathsurround=2pt
\pagewidth{6.5in}


\UseAMSsymbols
\input amssym.def
\loadbold


\def\etabar{{\bar{\eta}}}

\def\Fcal{{\Cal{F}}}
\def\Lcal{{\Cal{L}}}

\def\F{{\Bbb{F}}}
\def\Fbar{{\overline{\Bbb{F}}}}

\def\Qellbar{{\overline{\Bbb{Q}}_\ell}}
\def\Qellprime{{\Bbb{Q}_{\ell'}}}
\def\Qellprimebar{{\overline{\Bbb{Q}}_{\ell'}}}
\def\Z{{\Bbb{Z}}}
\def\Zhat{{\widehat{\Bbb{Z}}}}

\def\GG{{\operatorname{\bold{G}}}}
\def\Garith{{\operatorname{G_{arith}}}}
\def\Ggeom{{\operatorname{G_{geom}}}}

\def\charpoly{\operatorname{ch}}
\def\Ext{\operatorname{Ext}}
\def\Frob{\operatorname{Frob}}
\def\Gal{\operatorname{Gal}}
\def\GL{\operatorname{GL}}
\def\Hom{\operatorname{Hom}}
\def\K{\operatorname{K}}
\def\Ker{\operatorname{Ker}}

\def\rank{\operatorname{rk}}
\def\Tr{\operatorname{Tr}}
\def\W{\operatorname{W}}

\def\Res{\operatorname{Res}}
\def\Ind{\operatorname{Ind}}

\def\threeprime{{\mathsurround=0pt$3'$}}


\def\listitems{\begingroup \parskip=0pt \parindent=20pt}
\def\endlistitems{\par\endgroup}
\def\litem#1{\par \hangafter=1\hangindent=20pt \indent\llap{\rm #1\enspace}\ignorespaces}

\def\listiitems{\begingroup \parskip=0pt \parindent=40pt}
\def\endlistiitems{\par\endgroup}
\def\liitem#1{\par \hangafter=1\hangindent=40pt \indent\llap{\rm #1\enspace}\ignorespaces}

\def\qed{{\mathsurround=0pt\hfill$\bold\square$}}


\topmatter

\title 
Independence of $\boldsymbol{\ell}$ in Lafforgue's theorem
\endtitle

\rightheadtext{Independence of $\ell$ in Lafforgue's theorem}

\author
CheeWhye Chin
\endauthor


\address
Department of Mathematics, Princeton University,
Princeton, NJ 08544, U.S.A.
\endaddress

\curraddr
Department of Mathematics, University of California,
Berkeley, CA 94720, U.S.A.
\endcurraddr

\email 
cchin\@math.princeton.edu, cheewhye\@math.berkeley.edu
\endemail


\date
June 25, 2002
\enddate





\abstract
Let $X$ be a smooth curve
over a finite field of characteristic $p$,
let $\ell\ne p$ be a prime number,
and let $\Lcal$ be an irreducible lisse $\Qellbar$-sheaf on $X$
whose determinant is of finite order.
By a theorem of L.\,Lafforgue,
for each prime number $\ell'\ne p$,
there exists
an irreducible lisse $\Qellprimebar$-sheaf $\Lcal'$ on $X$
which is compatible with $\Lcal$,
in the sense that
at every closed point $x$ of $X$,
the characteristic polynomials of Frobenius at $x$
for $\Lcal$ and $\Lcal'$
are equal.
We prove
an ``independence of $\ell$'' assertion
on the fields of definition of
these irreducible $\ell'$-adic sheaves $\Lcal'$:
namely, that
there exists a number field $F$
such that
for any prime number $\ell'\ne p$,
the $\Qellprimebar$-sheaf $\Lcal'$ above
is defined over
the completion of $F$ at one of its $\ell'$-adic places.
\endabstract

%
%

\endtopmatter



\parskip=10pt plus 2pt



\document

\head{}
Introduction
\endhead

In the recent spectacular work \cite{L},
L.\,Lafforgue has proved
the Langlands Correspondence
and the Ramanujan-Petersson Conjecture
for $\GL_r$ over function fields.
As a consequence,
he has also established
the following fundamental result
concerning irreducible lisse $\ell$-adic sheaves
on curves over finite fields.

\proclaim{Theorem (L.~Lafforgue, {\rm \cite{L}~Th\'eor\`eme~VII.6})}
Let $X$ be a smooth curve
over a finite field of charac\-teristic~$p$.
Let $\ell \ne p$ be a prime number,
and let $\Lcal$ be a lisse $\Qellbar$-sheaf on $X$,
which is irreducible,
of rank~$r$,
and whose determinant is of finite order.
\listitems
\litem{(1)}
There exists a number field $E \subset \Qellbar$
such that
for every closed point $x$ of $X$,
the polynomial
$$ \det(1-T\,\Frob_x, \Lcal)
$$
has coefficients in $E$.
\litem{(2)}
Let $x$ be a closed point of $X$,
and let $\alpha\in\Qellbar$ be
an eigenvalue of Frobenius at~$x$
acting on $\Lcal$,
i.e.~$1/\alpha$ is a root of the polynomial
$$\det(1-T\,\Frob_x, \Lcal).
$$
Then:
\listiitems
\liitem{(a)}
$\alpha$ is an algebraic number;
\liitem{(b)}
for every archimedean absolute value $|\cdot|$ of $E(\alpha)$,
one has
$$|\alpha| = 1;
$$
\liitem{(c)}
for every non-archimedean valuation $\lambda$ of $E(\alpha)$
not lying over $p$,
$\alpha$ is a $\lambda$-adic unit,
i.e.~one has
$$ \lambda(\alpha) = 0;
$$
\liitem{(d)}
for every non-archimedean valuation $\nu$ of $E(\alpha)$
lying over $p$,
one has
$$ \left| \dfrac{\nu(\alpha)}{\nu(\#\kappa(x))} \right|
   \le
   \dfrac{(r-1)^2}{r}.
$$
\endlistiitems
\litem{(3)}
For any place $\lambda'$ of $E$
lying over a prime number $\ell'\ne p$,
and for any algebraic closure $\Qellprimebar$
of the completion $E_{\lambda'}$ of $E$ at $\lambda'$,
there exists
a lisse $\Qellprimebar$-sheaf $\Lcal'$ on $X$,
which is irreducible,
of rank $r$,
such that
for every closed point $x$ of $X$, one has
$$ \det(1-T\,\Frob_x, \Lcal')
 = \det(1-T\,\Frob_x, \Lcal)
	\qquad\text{(equality in $E[T]$)}.
$$
Moreover,
the sheaf $\Lcal'$
is defined over
a finite extension of $E_{\lambda'}$.
\endlistitems
\endproclaim

In part~(3) of Lafforgue's theorem,
it is not a~priori clear that
the number field $E$ may be
replaced by a finite extension (in $\Qellbar$)
so that
the various $\Qellprimebar$-sheaves $\Lcal'$
form an {\sl $(E,\Lambda)$-compatible system\/}
in the sense of Katz
(cf.~\cite{K}, pp.~202--203,
 ``The notion of $(E,\Lambda)$-compatibility''),
or equivalently, that
they form
an {\sl $E$-rational system of $\lambda$-adic representations\/}
in the sense of Serre
(cf.~\cite{Se}, \S2.3 and \S2.5).
The existence of a number field
with this property
may be interpreted as
an ``independence of $\ell$'' assertion on the
fields of definition of
these irreducible $\ell'$-adic sheaves $\Lcal'$.
We shall prove that
this is indeed the case.

\proclaim{Theorem}
With the notation and hypotheses
of Lafforgue's Theorem,
the following assertion holds.
\listitems
\litem{(\threeprime)}
There exists
a finite extension $F$ of $E$ in $\Qellbar$
such that
for any place $\lambda'$ of the number field $F$
lying over a prime number $\ell'\ne p$,
there exists
a lisse $F_{\lambda'}$-sheaf $\Lcal'$ on $X$
(i.e. a lisse $\Qellprimebar$-sheaf defined over $F_{\lambda'}$),
which is absolutely irreducible,
of rank $r$,
such that
for every closed point $x$ of $X$, one has
$$ \det(1-T\,\Frob_x, \Lcal')
 = \det(1-T\,\Frob_x, \Lcal)
	\qquad\text{(equality in $E[T]$)}.
$$
\endlistitems
\endproclaim

According to a conjecture of Deligne
(cf. \cite{D} Conjecture~(1.2.10)),
all four assertions (1),(2),(3),(\threeprime)
should also hold
in the general case
when $X$ is a normal variety of arbitrary dimension
over a finite field.
Our proof of assertion~(\threeprime)
uses assertions~(1) and~(3) of Lafforgue's Theorem
only as ``black boxes'';
so assertion~(\threeprime) will hold for 
higher dimensional varieties
{\it if\/} parts~(1) and~(3) of Lafforgue's theorem
hold for these varieties.
To state this more precisely,
we make assertions~(1) and~(3) into hypotheses,
as follow.

\proclaim{Definition}{\rm
Let $\F_q$ be a finite field of characteristic $p$,
and let $\ell\ne p$ be a prime number.
Let $Y$ be a normal variety over $\F_q$,
and let $\Fcal$ be a lisse $\Qellbar$-sheaf on $Y$,
which is irreducible,
and whose determinant is of finite order.
We shall say that
{\sl hypothesis~(1) holds for $(Y,\Fcal)$\/}
if:
\listitems
\litem{(1)}
there exists a number field $E \subset \Qellbar$
such that
for every closed point $y$ of $Y$,
the polynomial
$$ \det(1-T\,\Frob_y, \Fcal)
$$
has coefficients in $E$.
\endlistitems
\noindent
When hypothesis~(1) holds for $(Y,\Fcal)$,
we shall say that
{\sl hypothesis~(3) holds for $(Y,\Fcal)$\/}
if:
\listitems
\litem{(3)}
for any place $\lambda'$ of $E$
lying over a prime number $\ell'\ne p$,
and for any algebraic closure $\Qellprimebar$
of the completion $E_{\lambda'}$ of $E$ at $\lambda'$,
there exists
a lisse $\Qellprimebar$-sheaf $\Fcal'$ on $Y$,
which is irreducible,
such that
for every closed point $y$ of $Y$, one has
$$ \det(1-T\,\Frob_y, \Fcal')
 = \det(1-T\,\Frob_y, \Fcal)
	\qquad\text{(equality in $E[T]$)}.
$$
\endlistitems
}
\endproclaim

With this definition,
our goal is to prove:

\proclaim{Main Theorem}
Let $\F_q$ be a finite field of characteristic $p$,
and let $\ell\ne p$ be a prime number.
Let $X$ be a normal variety over $\F_q$.
Assume that:
\block
   for any normal variety $Y$ over $\F_q$
   which is finite etale over $X$,
   and for any lisse $\Qellbar$-sheaf $\Fcal$ on $Y$,
   which is irreducible,
   and whose determinant is of finite order,
   hypotheses~(1) and~(3) hold for the pair $(Y,\Fcal)$.
\endblock
Let $\Lcal$ be a lisse $\Qellbar$-sheaf on $X$,
which is irreducible,
of rank~$r$,
and whose determinant is of finite order.
Let $E\subset\Qellbar$ denote
the number field given by hypothesis~(1) applied to $(X,\Lcal)$.
Then:
\listitems
\litem{(\threeprime)}
There exists
a finite extension $F$ of $E$ in $\Qellbar$
such that
for any place $\lambda'$ of the number field $F$
lying over a prime number $\ell'\ne p$,
there exists
a lisse $F_{\lambda'}$-sheaf $\Lcal_{\lambda'}$ on $X$,
which is absolutely irreducible,
of rank $r$,
such that
for every closed point $x$ of $X$, one has
$$ \det(1-T\,\Frob_x, \Lcal_{\lambda'})
 = \det(1-T\,\Frob_x, \Lcal)
	\qquad\text{(equality in $E[T]$)}.
$$
\endlistitems
\endproclaim

We shall prove this theorem
by exploiting properties of
the {\it monodromy groups\/} associated to
these irreducible lisse sheaves.
The proof begins in \S4,
after a discussion of the preliminary results we need:
propo\-sitions~1 and~2 of \S1,
corollary~6 of \S2,
and propositions~7 and~9 of \S3.

I am grateful to P.\,Deligne, J.\,de\,Jong and especially N.\,Katz
for many fruitful conversations,
from which I have learned much of the material presented here.

\head{\S1.}
Monodromy Groups
\endhead

In this section, we recall
some basic properties
of monodromy groups of lisse $\ell$-adic sheaves
on varieties over a finite field;
see \cite{D}~\S1.1 and \S1.3 for details.

Let $X$ be a normal, geometrically connected variety
over a finite field $\F_q$ of characteristic~$p$.
Let $\etabar @>>> X$ be a geometric point of $X$,
and let $\Fbar_q$ be
the algebraic closure $\F_q$ in $\kappa(\etabar)$;
we regard $\etabar$ also as a geometric point
of $X\otimes_{\F_q}\Fbar_q$.
The profinite groups
$\pi_1(X,\etabar)$ and $\pi_1(X\otimes_{\F_q}\Fbar_q,\etabar)$
are respectively called
the {\sl arithmetic fundamental group\/} of $X$
and the {\sl geometric fundamental group\/} of $X$.
They sit in a short exact sequence
$$ 1 @>>>
   \pi_1(X\otimes_{\F_q}\Fbar_q,\etabar) @>>>
   \pi_1(X,\etabar) @>{\deg}>>
   \Gal(\Fbar_q/\F_q) @>>> 1.
$$
The group $\Gal(\Fbar_q/\F_q)$
has a canonical topological generator $\Frob_{\F_q}$ called
the {\sl geometric Frobenius\/},
which is defined as the inverse of
the {\sl arithmetic Frobenius\/} automorphism
$ a \mapsto a^q $ of the field $\Fbar_q$.
We have the canonical isomorphism
$$ \Zhat @>{\cong}>> \Gal(\Fbar_q/\F_q),
   \qquad\text{sending $1$ to $\Frob_{\F_q}$}.
$$

For a prime number $\ell\ne p$, the functor
$$\matrix\format\r&\ \c\ &\l			\\
 \{ \text{lisse $\Qellbar$-sheaves on $X$} \}
 & \longrightarrow &	
 \{ \text{finite dimensional continuous
	  $\Qellbar$-representations
	  of $\pi_1(X,\etabar)$} \}	\\
 \Lcal \qquad
 & \mapsto &
 \qquad \Lcal_\etabar				\\
\endmatrix
$$
is an equivalence of categories;
a similar statement holds with $X\otimes_{\F_q}\Fbar_q$
in place of $X$.
Via this equivalence,
standard notions associated to representations
(e.g.~irreducibility, semisimplicity, constituent, etc.)
are also applicable to lisse sheaves.

Let $\Lcal$ be a lisse $\Qellbar$-sheaf on $X$,
corresponding to the continuous monodromy representation
$$ \pi_1(X,\etabar) @>>> \GL(\Lcal_\etabar)
$$
of the arithmetic fundamental group of $X$.
The {\sl arithmetic monodromy group\/}
$\Garith(\Lcal,\etabar)$ of $\Lcal$
is the Zariski closure of the image of
$\pi_1(X,\etabar)$ in $\GL(\Lcal_\etabar)$.
The inverse image $\Lcal\otimes_{\F_q}\Fbar_q$
of $\Lcal$ on $X\otimes_{\F_q}\Fbar_q$
is a lisse $\Qellbar$-sheaf on $X\otimes_{\F_q}\Fbar_q$,
corresponding to the continuous monodromy representation
$$ \pi_1(X\otimes_{\F_q}\Fbar_q,\etabar)
	\hookrightarrow \pi_1(X,\etabar) @>>> \GL(\Lcal_\etabar)
$$
of the geometric fundamental group of $X$,
obtained by restriction.
The {\sl geometric monodromy group\/}
$\Ggeom(\Lcal,\etabar)$ of $\Lcal$
is the Zariski closure of the image of
$\pi_1(X\otimes_{\F_q}\Fbar_q,\etabar)$ in $\GL(\Lcal_\etabar)$.

Both $\Garith(\Lcal,\etabar)$ and $\Ggeom(\Lcal,\etabar)$
are linear algebraic groups,
and it is clear that
$\Ggeom(\Lcal,\etabar)$ is
a closed normal subgroup of $\Garith(\Lcal,\etabar)$.
Both $\Garith(\Lcal,\etabar)$ and $\Ggeom(\Lcal,\etabar)$
are given with a faithful represen\-tation on $\Lcal_\etabar$
corresponding to their realizations
as subgroups of $\GL(\Lcal_\etabar)$.
Thus, if $\Lcal$ is semisimple
(as a representation of $\pi_1(X,\etabar)$,
 and therefore as a representation of
 $\pi_1(X\otimes_{\F_q}\Fbar_q,\etabar)$),
then both $\Garith(\Lcal,\etabar)$ and $\Ggeom(\Lcal,\etabar)$
are (possibly non-connected) reductive algebraic groups.

\proclaim{Proposition~1}
Let $\Lcal$ be a lisse $\Qellbar$-sheaf on $X$.
\listitems
\litem{(i)}
If $\Lcal$ is semisimple,
then $\Ggeom(\Lcal,\etabar)$ is a (possibly non-connected)
semisimple algebraic group.
\litem{(ii)}
If $\Lcal$ is irreducible,
and its determinant is of finite order,
then $\Garith(\Lcal,\etabar)$ is a (possibly non-connected)
semisimple algebraic group,
containing $\Ggeom(\Lcal,\etabar)$ as a normal subgroup
of finite index.
\endlistitems
\endproclaim

Assertion~(i) is \cite{D}~Corollaire~(1.3.9).
For the proof of assertion~(ii),
we shall make use of the construction in \cite{D}~(1.3.7),
which we summarize below.

Recall that
the {\sl Weil group\/} $\W(\Fbar_q/\F_q)$ of $\F_q$
is the subgroup of $\Gal(\Fbar_q/\F_q)$ consisting of
integer-powers of $\Frob_{\F_q}$;
it is considered as a topological group given with
the discrete topology,
and we have the canonical isomorphism
$$ \Z @>{\cong}>> \W(\Fbar_q/\F_q),
	\qquad\text{sending $1$ to $\Frob_{\F_q}$}.
$$
The {\sl Weil group\/} $\W(X,\etabar)$ of $X$
is the preimage of $\W(\Fbar_q/\F_q)$ in $\pi_1(X,\etabar)$ 
by the degree homomorphism
$ \pi_1(X,\etabar) @>{\deg}>> \Gal(\Fbar_q/\F_q) $;
it is considered as a topological group given with
the product topology via the isomorphism
$$ \W(X,\etabar) \cong
   \pi_1(X\otimes_{\F_q}\Fbar_q,\etabar)
        \rtimes_{\Gal(\Fbar_q/\F_q)} \W(\Fbar_q/\F_q),
$$
where $\pi_1(X\otimes_{\F_q}\Fbar_q,\etabar)$
retains its profinite topology,
and is an open and closed subgroup of $\W(X,\etabar)$.
These groups sit in the following diagram:
$$\CD
   1 @>>> \pi_1(X\otimes_{\F_q}\Fbar_q,\etabar)
     @>>> \W(X,\etabar)
     @>{\deg}>> \Z \cong \W(\Fbar_q/\F_q)
     @>>> 1					\\
   @.	  \vert\vert@.
	  @V{\cap}VV
	  @V{\cap}VV				\\
   1 @>>> \pi_1(X\otimes_{\F_q}\Fbar_q,\etabar)
     @>>> \pi_1(X,\etabar)
     @>{\deg}>> \Zhat \cong \Gal(\Fbar_q/\F_q)
     @>>> 1
\endCD
$$
where the right two vertical arrows are
inclusion homomorphisms with dense images.
(Note that
 the topo\-logies of $\W(X,\etabar)$ and $\W(\Fbar_q/\F_q)$
 are not the ones induced by
 the right two vertical arrows!)

Given a lisse $\Qellbar$-sheaf $\Lcal$ on $X$,
the {\it push-out construction\/} of \cite{D}~(1.3.7) produces
an algebraic group $\GG(\Lcal,\etabar)$,
which is locally of finite type, but not quasi-compact;
it is characterized by the fact that
it sits in a diagram:
$$\CD
   1 @>>> \pi_1(X\otimes_{\F_q}\Fbar_q,\etabar)
     @>>> \W(X,\etabar)
     @>{\deg}>> \Z \cong \W(\Fbar_q/\F_q)
     @>>> 1					\\
   @.	  @VVV
	  @VVV
	  \vert\vert@.				\\
   1 @>>> \Ggeom(\Lcal,\etabar)
     @>>> \GG(\Lcal,\etabar)
     @>{\deg}>> \Z \cong \W(\Fbar_q/\F_q)
     @>>> 1					\\
   @.	  @.
	  @VVV					\\
     @.
     @.   \GL(\Lcal_\etabar)
\endCD
$$
such that the composite of the two continuous homomorphisms
$$ \W(X,\etabar) @>>> \GG(\Lcal,\etabar) @>>> \GL(\Lcal_\etabar)
$$
is equal to
the continuous representation of $\W(X,\etabar)$ on $\Lcal_\etabar$
obtained via restriction:
$$ \W(X,\etabar) \hookrightarrow \pi_1(X,\etabar) @>>> \GL(\Lcal_\etabar).
$$

\demo{Proof of Proposition~1~(ii)}
From assertion~(i), we already know that
the group $\Ggeom(\Lcal,\etabar)$
is a semisimple closed normal subgroup
of $\Garith(\Lcal,\etabar)$.
Hence, to prove assertion~(ii),
it suffices for us to show that
$\Garith(\Lcal,\etabar)$ contains $\Ggeom(\Lcal,\etabar)$
as a subgroup of finite index,
for then both groups will have the same identity component,
which is a connected semisimple algebraic group.

Since $\W(X,\etabar) \hookrightarrow \pi_1(X,\etabar)$
is an inclusion with dense image,
$\Garith(\Lcal,\etabar)$ can also be described as
the Zariski closure of
the image of $\W(X,\etabar)$ in $\GL(\Lcal_\etabar)$;
likewise, since $\W(X,\etabar) \hookrightarrow \GG(\Lcal,\etabar)$
is an inclusion with dense image,
$\Garith(\Lcal,\etabar)$ is also equal to
the Zariski closure of
the image of $\GG(\Lcal,\etabar)$ in $\GL(\Lcal_\etabar)$.
Let
$$ \rho: \GG(\Lcal,\etabar) @>>> \GL(\Lcal_\etabar)
$$
denote the canonical homomorphism
from $\GG(\Lcal,\etabar)$ into $\GL(\Lcal_\etabar)$;
then the composite map
$$ \Ggeom(\Lcal,\etabar)
     \hookrightarrow \GG(\Lcal,\etabar) @>{\rho}>> \GL(\Lcal_\etabar)
$$
is just the identity map on $\Ggeom(\Lcal,\etabar)$.
We are thus reduced to showing that
$\rho^{-1}(\Ggeom(\Lcal,\etabar))$ is
a subgroup of $\GG(\Lcal,\etabar)$ of finite index.

The fundamental fact we need about $\GG(\Lcal,\etabar)$
is \cite{D}~Corollaire~(1.3.11), 
which asserts that
because $\Lcal$ is irreducible (hence semisimple) by hypothesis,
there exists some element $g$
in the center of $\GG(\Lcal,\etabar)$
whose degree is $>0$
(i.e. $g$ maps to a positive integer under
 $\GG(\Lcal,\etabar) @>{\deg}>> \Z \cong \W(\Fbar_q/\F_q)$).
Therefore, $\rho(g)$ is an element of $\GL(\Lcal_\etabar)$
which centralizes $\rho(\GG(\Lcal,\etabar))$,
and so it centralizes $\Garith(\Lcal,\etabar)$.
Since $\Lcal$ is irreducible
as a representation of $\pi_1(X,\etabar)$
and hence
as a representation of $\Garith(\Lcal,\etabar)$,
it follows that
$\rho(g)$ must be a scalar.

By hypothesis,
the determinant of $\Lcal$ is of finite order,
which means that
the $1$-dimensional representation
of $\pi_1(X,\etabar)$ on
the determinant $\det(\Lcal_\etabar)$ of $\Lcal_\etabar$
is given by a character of finite order, say $d$.
The same is therefore true for $\det(\Lcal_\etabar)$ as
a representation of $\W(X,\etabar)$ and of $\GG(\Lcal,\etabar)$.
From this it follows that,
if $\Lcal$ has rank~$r$,
then $\rho(g)$ is a scalar which is
a root of unity of order dividing $d\,r$,
and so $g^{\,d\,r} \in \GG(\Lcal,\etabar)$ lies in the kernel of $\rho$.
Hence $\rho^{-1}(\Ggeom(\Lcal,\etabar))$ contains
$\deg^{-1}(\deg(g^{\,d\,r}))$ in $\GG(\Lcal,\etabar)$,
which is of finite index in $\GG(\Lcal,\etabar)$.
\qed
\enddemo

Let $\Lcal$ be a lisse $\Qellbar$-sheaf $\Lcal$ on $X$.
Its arithmetic monodromy group $\Garith(\Lcal,\etabar)$
contains the identity component $\Garith(\Lcal,\etabar)^0$
as an open normal subgroup;
$\Garith(\Lcal,\etabar)^0$ is a connected algebraic group.
The faithful representation
$$ \Garith(\Lcal,\etabar) \hookrightarrow \GL(\Lcal_\etabar)
$$
of $\Garith(\Lcal,\etabar)$,
when restricted to
the subgroup $\Garith(\Lcal,\etabar)^0$ of $\Garith(\Lcal,\etabar)$,
gives a faithful representation
$$ \Garith(\Lcal,\etabar)^0 \hookrightarrow
   \Garith(\Lcal,\etabar) \hookrightarrow \GL(\Lcal_\etabar)
$$
of $\Garith(\Lcal,\etabar)^0$ on $\Lcal_\etabar$.
We say that
the lisse sheaf $\Lcal$ is {\sl Lie-irreducible\/}
if $\Lcal_\etabar$ is irreducible
as a representation of $\Garith(\Lcal,\etabar)^0$.
It is clear that
Lie-irreducibility implies irreducibility.

\proclaim{Proposition~2}
Let $\Lcal$ be a lisse $\Qellbar$-sheaf on $X$,
which is Lie-irreducible,
and whose determinant is of finite order.
Then there exist
$\alpha\in\Qellbar$
and a closed point $x_0$ of $X$,
such that
$\alpha$ is an eigenvalue of multiplicity~one
of $\Frob_{x_0}$ acting on $\Lcal$;
i.e.
$1/\alpha$ is a root of multiplicity~one
of the polynomial
$$ \det(1-T\,\Frob_{x_0},\Lcal).
$$
\endproclaim

\demo{Proof of Proposition~2}
First, we claim that
it is a Zariski-open condition
for an element of $\Garith(\Lcal,\etabar)$
to have an eigenvalue of multiplicity~one
on $\Lcal_\etabar$;
in other words, we claim that
the set
$$ U := \{g \in\Garith(\Lcal,\etabar) \ :\ 
 	    \text{$g$ acting on $\Lcal_\etabar$
		  has an eigenvalue
		  of multiplicity~one in $\Qellbar$} \}
$$
is a {\it Zariski-open\/} subset of $\Garith(\Lcal,\etabar)$.
We show this as follows.
For an element $g \in\Garith(\Lcal,\etabar)$,
let $\charpoly(g)\in\Qellbar[T]$ denote
the characteristic polynomial of $g$;
then the set $U$ can also be described as
$$U =   \{ g \in\Garith(\Lcal,\etabar) \ :\ 
	  \text{$\charpoly(g) \in\Qellbar[T]$
		has a root of multiplicity~one
		in $\Qellbar$} \}.
$$
Let $r$ be the rank of $\Lcal_\etabar$;
then $\charpoly$ gives rise to a morphism of $\Qellbar$-varieties
$$ \charpoly: \Garith(\Lcal,\etabar) @>>> \Qellbar[T]^{\text{monic}}_{\deg r},
	\qquad g \mapsto \charpoly(g),
$$
where $\Qellbar[T]^{\text{monic}}_{\deg r}$ denotes
the affine space of monic polynomials in $T$ of degree~$r$.
For $g \in\Garith(\Lcal,\etabar)$,
the polynomial $\charpoly(g)$
has a root of multiplicity~one in $\Qellbar$
if and only if
it does not divide
the square $\charpoly(g)^{\prime\,2}$
of its derivative $\charpoly(g)'$ in $\Qellbar[T]$.
Thus it suffices for us to show that
the set
$$ Z := \{ f \in \Qellbar[T]^{\text{monic}}_{\deg r} \ :\ 
		\text{$f$ divides $f^{\prime\,2}$ in $\Qellbar[T]$} \}
$$
is Zariski-closed in $\Qellbar[T]^{\text{monic}}_{\deg r}$.
But for $f \in \Qellbar[T]^{\text{monic}}_{\deg r}$,
the Euclidean division algorithm shows that
the remainder of dividing $f^{\prime\,2}$ by $f$
is a polynomial of degree~$< r$
whose coefficients are given by
certain (universal) $\Z$-polynomial expressions
in terms of the coefficients of $f$;
as the set $Z$ above
is precisely the zero-set of these polynomial expressions,
it is Zariski-closed.

Next, we claim that
the set $U$ above
is in fact Zariski-open and {\it non-empty\/}
in $\Garith(\Lcal,\etabar)$.
Indeed,
by part~(ii) of proposition~1,
$\Garith(\Lcal,\etabar)^0$ is
a connected semisimple algebraic group;
the representation $\Lcal_\etabar$
of $\Garith(\Lcal,\etabar)^0$
is irreducible by hypothesis,
and so
by the representation theory of connected semisimple algebraic groups,
it is classified by
its highest weight,
which occurs with multiplicity~one.
Thus a generic element of
any maximal torus of $\Garith(\Lcal,\etabar)^0$
lies in $U$.

Finally, by \u{C}ebotarev's density theorem,
there exist infinitely many closed points $x$ of $X$
whose Frobenius conjugacy classes $\Frob_x \subset \pi_1(X,\etabar)$
are mapped into $U$ under
the monodromy representation of $\pi_1(X,\etabar)$ on $\Lcal_\etabar$.
Thus we can
pick $x_0$ to be any one of these closed points of $X$,
and pick $\alpha\in\Qellbar$ to be
an eigenvalue of multiplicity~one
of $\Frob_{x_0}$ acting on $\Lcal$.
\qed
\enddemo

\remark{Remark}
In proposition~2,
it is not enough to just assume that
the lisse $\Qellbar$-sheaf $\Lcal$ is irreducible;
the assumption that it is {\it Lie-irreducible\/}
is necessary.
If $\Lcal$ is irreducible but not Lie-irreducible,
it may happen that
every element of $\Garith(\Lcal,\etabar)$
acting on $\Lcal_\etabar$
has repeated eigenvalues,
which is to say that
the set $U \subset \Garith(\Lcal,\etabar)$ 
in the proof of the proposition
is empty.
For a specific example,
we may take $\Garith(\Lcal,\etabar)$ to be
the finite symmetric group on $6$~letters,
and take $\Lcal_\etabar$ to be
the $16$-dimensional irreducible representation
of this finite group;
such a situation can arise geometrically.
\endremark

\head{\S2.}
D\'evissage of Representations
\endhead

Let $k$ be an algebraically closed field
of characteristic~$0$ ---
such as $\Qellbar$.
In this section, we consider
(possibly non-connected) reductive groups over $k$
and their finite dimensional
$k$-rational representations.
If $G$ is such a reductive group,
any $k$-rational representation of $G$ is semisimple
(a direct sum of irreducible representations),
since $k$ is of characteristic~$0$.
By the quasi-compactness of $G$,
a subgroup $H$ of $G$ is (Zariski-) open
if and only if it is (Zariski-) closed of finite index,
in which case
$H$ necessarily contains the identity component $G^0$ of $G$.

The following two results are proved in \cite{I}
for representations of finite groups.
The same proofs, with minor modifications,
work for representations of reductive groups.
We reproduce the (modified) arguments below
for the sake of completeness.

\proclaim{Lemma~3 (I.~M.~Isaacs, {\rm \cite{I}~Theorem~6.18})}
Let $G$ be a reductive group,
and let $K$ and $L$ be open normal subgroups of $G$,
with $L \subseteq K$.
Suppose that $K/L$ is abelian,
and that there does not exist a normal subgroup $M$ of $G$
with $L \subsetneq M \subsetneq K$.
Let $\pi$ be an irreducible representation of $K$
whose isomorphism class is
invariant under $G$-conjugation.
Then one of the following holds:
\listitems
\litem{(i)}
$\Res^K_L(\pi)$ is isomorphic to a direct sum
$\sigma_1 \oplus\cdots\oplus \sigma_t$
of $t := [K:L]$ many irreducible representations
$\sigma_1,\ldots,\sigma_t$ of $L$
which are pairwise non-isomorphic;
\litem{(ii)}
$\Res^K_L(\pi)$ is an irreducible representation of $L$;
\litem{(iii)}
$\Res^K_L(\pi)$ is isomorphic to $\sigma^{\oplus e}$,
where $\sigma$ is an irreducible representation of $L$,
and $e^2 = [K:L]$.
\endlistitems
\endproclaim

\demo{Proof of Lemma~3}
Since $L$ is normal in $K$,
the irreducible constituents of $\Res^K_L(\pi)$
are $K$-conjugate to one another,
and each of these constituents occurs in $\Res^K_L(\pi)$
with the same multiplicity.
Choose any irreducible constituent $\sigma$ of $\Res^K_L(\pi)$,
and let
$$ I := \{ g\in G \ :\ {}^g\sigma \cong \sigma
		\quad\text{as representations of $L$} \}
$$
be the open subgroup of $G$ (containing $L$)
which stabilizes the isomorphism type of $\sigma$
under $G$-conjugation.
Since $\pi$ is invariant under $G$-conjugation,
every $G$-conjugate of $\sigma$
is a constituent of $\Res^K_L(\pi)$,
and so every $G$-conjugate of $\sigma$
is $K$-conjugate to $\sigma$.
It follows that $[G:I] = [K:K\cap I]$,
and hence $KI = G$.
Since $K/L$ is abelian,
$K\cap I$ is normal in $K$;
since $K$ is normal in $G$,
$K\cap I$ is normal in $I$.
As $KI = G$,
we see that
$K\cap I$ is normal in $G$.
From the hypothesis of the proposition,
it follows that
$K \cap I$ is either $L$ or $K$.

Suppose $K \cap I = L$.
Then there are $t = [K:L]$ many
pairwise non-isomorphic irreducible constituents
$\sigma = \sigma_1, \ldots, \sigma_t$ of $\Res^K_L(\pi)$,
and so we have
$$ \Res^K_L(\pi) \cong
 (\sigma_1 \oplus\cdots\oplus \sigma_t)^{\oplus e}
$$
for some multiplicity $e\ge 1$.
The constituents $\sigma_j$ of $\Res^K_L(\pi)$
are $K$-conjugate to one another,
and so they have the same rank as $\sigma$.
Hence
$$ \rank(\pi) = \rank(\Res^K_L(\pi)) = e\,t\,\rank(\sigma).
$$
But $\pi$ is a constituent of $\Ind^K_L(\sigma)$,
so
$$ \rank(\pi) \le \rank(\Ind^K_L(\sigma)) = t\,\rank(\sigma).
$$
Thus $e=1$, and this is case~(i).

Henceforth suppose $K \cap I = K$.
Then $\sigma$ is invariant under $K$-conjugation,
so we have
$$ \Res^K_L(\pi) \cong \sigma^{\oplus e}
$$
for some multiplicity $e\ge 1$.
Let $\chi_1,\ldots,\chi_t$ be the distinct linear
characters of the abelian group $K/L$.
Then $\chi_1\otimes\pi,\ldots,\chi_t\otimes\pi$ are
irreducible representations of $K$,
each having the same rank as $\pi$,
and we have
$$ \Res^K_L(\chi_j\otimes\pi) \cong \sigma^{\oplus e}
		\qquad\text{for each $j=1,\ldots,t$}.
$$

Suppose $\chi_1\otimes\pi,\ldots,\chi_t\otimes\pi$
are pairwise non-isomorphic representations of $K$.
Then we obtain an inclusion
$$ \bigoplus_{j=1}^t (\chi_j\otimes\pi)^{\oplus e}
	\subseteq \Ind^K_L(\sigma).
$$
Comparing ranks, we get
$$ e\,t\,\rank(\pi) \le \rank(\Ind^K_L(\sigma)) = t\,\rank(\sigma),
$$
and so
$$ e\,\rank(\pi) \le \rank(\sigma).
$$
But
$$ e\,\rank(\sigma)
 = \rank(\Res^K_L(\pi))
 = \rank(\pi).
$$
Thus $e=1$, and this is case~(ii).

In the remaining situation,
at least two of the representations
$\chi_1\otimes\pi,\ldots,\chi_t\otimes\pi$
are isomorphic;
this implies that $\pi \cong \chi\otimes\pi$
for some non-trivial linear character $\chi$ of $K/L$.
Let $M = \Ker(\chi)$;
we have $L \subseteq M \subsetneq K$.
First, consider the representation $\pi$,
with trace-function
$$ \Tr\circ\pi : K @>>> k,
	\qquad x \mapsto \Tr(\pi(x)).
$$
On $K-M$, the linear character $\chi$ takes values different from~$1$;
since
$\Tr\circ\pi
 = \Tr\circ(\chi\otimes\pi)
 = \chi\cdot(\Tr\circ\pi)$,
it follows that $\Tr\circ\pi$ vanishes on $K-M$.
Since the representation $\pi$
is invariant under $G$-conjugation,
it follows that $\Tr\circ\pi$ vanishes on $K-gMg^{-1}$
for all $g\in G$.
The normal subgroup
$\bigcap_{g\in G} gMg^{-1}$ of $G$
contains $L$ and is properly contained in $K$,
so it must be equal to $L$ by hypothesis.
Thus $\Tr\circ\pi$ vanishes on $K-L$.
Next, consider the representation 
$\Ind^K_L(\Res^K_L(\pi)) \cong \Ind^K_L(1)\otimes\pi$,
with its trace-function
$$ \Tr\circ\Ind^K_L(\Res^K_L(\pi)) : K @>>> k,
	\qquad x \mapsto \Tr(\Ind^K_L(1)(x))\,\Tr(\pi(x)).
$$
Since the trace-function of $\Ind^K_L(1)$
is $0$ on $K-L$ and is $t$ on $L$,
it follows that
the trace-function of $\Ind^K_L(\Res^K_L(\pi))$
vanishes on $K-L$,
and its values on $L$
are $t$ times those of $\Tr\circ\pi$.
Comparing the trace-functions of $\pi$ and $\Ind^K_L(\Res^K_L(\pi))$,
we see that
$$ \Tr\circ(\pi^{\oplus t}) = \Tr\circ\Ind^K_L(\Res^K_L(\pi)).
$$
By the trace comparison theorem of Bourbaki
(cf.~\cite{B}~\S12, no.~1, Prop.~3),
this implies
$$ \pi^{\oplus t} \cong \Ind^K_L(\Res^K_L(\pi))
$$
as representations of $K$.
Hence
$$ e^2 = \dim \Hom_L(\Res^K_L(\pi),\Res^K_L(\pi))
       = \dim \Hom_K(\pi,\Ind^K_L(\Res^K_L(\pi)))
       = t = [K:L]
$$
and this is case~(iii).
\qed
\enddemo

\proclaim{Proposition~4 (I.~M.~Isaacs, {\rm \cite{I}~Theorem~6.22})}
Let $G$ be a reductive group,
and let $N$ be an open normal subgroup of $G$
such that $G/N$ is a nilpotent finite group.
Let $\rho$ be an irreducible representation of $G$.
Then there exists an open subgroup $H$ of $G$
with $N \subseteq H \subseteq G$,
and an irreducible representation $\sigma$ of $H$,
such that
$\rho \cong \Ind^G_H(\sigma)$,
and such that
$\Res^H_N(\sigma)$ is
an irreducible representation of $N$.
\endproclaim

\remark{Remark}
The proposition holds in slightly greater generality:
we need only to assume that
$G/N$ is a solvable finite group
whose {\it chief factors\/} are of square-free orders;
see \cite{I}.
This technical condition is automatically verified
when $G/N$ is nilpotent or supersolvable.
\endremark

\demo{Proof of Proposition~4}
The theorem is clear when $G=N$.
We proceed by induction on $\#(G/N)$;
hence assume that the theorem holds
for any proper subgroup of $G$ containing $N$.
If $\Res^G_N(\rho)$ is irreducible,
then the theorem holds with $H = G$ and $\sigma = \rho$.
Hence suppose $\Res^G_N(\rho)$ is reducible.

Since $G/N$ is finite,
we can find an open normal subgroup $K$ of $G$
which is minimal for the conditions that
$N \subseteq K$ and $\Res^G_K(\rho)$ is irreducible.
Then $N \subsetneq K$ necessarily,
and so we can find an open normal subgroup $L$ of $G$
which is maximal for the conditions that
$N \subseteq L \subsetneq K$.
Since $G/N$ is nilpotent,
it follows that $K/L$ is cyclic of prime order, say $t$.

The isomorphism class of the irreducible representation
$\pi = \Res^G_K(\rho)$ of $K$
is invariant under $G$-conjugation,
since $\pi$ is the restriction of
an irreducible representation $\rho$ of $G$.
Thus we may apply lemma~3 to
the representation $\pi$ of $K$.
By the choice of $L$ and $K$,
$\Res^K_L(\pi)$ is not irreducible,
so case~(ii) cannot occur;
since $t = [K:L]$ is a prime number,
case~(iii) cannot occur.
Hence we are in case~(i), and it follows that
$\Res^G_L(\rho)$ is isomorphic to a direct sum
$\sigma_1 \oplus\cdots\oplus \sigma_t$ of $t$ many
irreducible representations $\sigma_1,\ldots,\sigma_t$
of $L$ which are pairwise non-isomorphic.

Let
$$ I := \{ g\in G \ :\ {}^g\sigma_1 \cong \sigma_1
		\quad\text{as representations of $L$} \}
$$
be the open subgroup of $G$ (containing $L$)
which stabilizes the isomorphism type of $\sigma_1$
under $G$-conjugation.
Thus $[G:I] = t$ is $>1$,
and $\rho \cong \Ind^G_I(\rho')$
for some irreducible representation $\rho'$ of $I$.
Applying the induction hypothesis to $I$,
we obtain an open subgroup $H$ of $I$
with $N \subseteq H \subseteq I$,
and an irreducible representation $\sigma$ of $H$,
such that $\rho' \cong \Ind^I_H(\sigma)$
and $\Res^H_N(\sigma)$ is
an irreducible representation of $N$.
Then $\rho \cong \Ind^G_H(\sigma)$,
which completes the proof of the proposition.
\qed
\enddemo

If $G$ is a reductive group over $k$,
we let $\K(G)$ denote
the Grothendieck group of the abelian category of
finite dimensional $k$-rational representations
of $G$.
It is clear that $\K(G)$ as a $\Z$-module
is freely generated by
the irreducible representations of $G$.
%
%
The tensor product of representations
gives rise to a commutative ring structure on $\K(G)$,
whose unit element is the class $1$
of the trivial representation of $G$.
If $H\subseteq G$ is an open subgroup,
then induction of representations from $H$ to $G$
gives rise to a homomorphism of $\Z$-modules
$$	\Ind : \K(H) @>>> \K(G).
$$
The projection formula shows that
the $\Ind$-image of $\K(H)$ in $\K(G)$
is an ideal.

Recall that, for $p$ a prime number,
a finite group $G$ is called {\sl $p$-elementary\/}
if it is isomorphic to a direct product $A\times B$,
where $A$ is a cyclic group of order prime to $p$,
and $B$ is a $p$-group.
A finite group $G$ is called {\sl elementary\/}
if it is $p$-elementary for some prime number~$p$.
It is clear that an elementary finite group
is nilpotent.

Let $G$ be a reductive group,
and $N$ be an open normal subgroup of $G$.
We say that, for a prime number~$p$,
an open subgroup $H$ of $G$
is {\sl $p$-elementary modulo $N$\/}
if one has the inclusions
$N \subseteq H \subseteq G$
and furthermore the finite quotient $H/N$
is $p$-elementary;
we say that $H$ is {\sl elementary modulo $N$\/}
if it is $p$-elementary modulo $N$
for some prime number~$p$.

\proclaim{Proposition~5 (R.~Brauer)}
Let $G$ be a reductive group,
and let $N$ be an open normal subgroup of $G$.
Then the $\Z$-homomorphism
$$ \Ind: \bigoplus_{\overset{H\subseteq G}\to{\text{elem.mod $N$}}}
	  \K(H)
	 @>>> \K(G)
$$
is surjective
(the direct sum is over all subgroups $H$ of $G$
 which are elementary modulo $N$).
\endproclaim

\demo{Proof of Proposition~5}
Recall that
Brauer's theorem on induced characters for finite groups
(see \cite{I}~Theorem~8.4 or \cite{H}~Theorem~34.2 for instance)
states that
if $G$ is a finite group,
then the $\Z$-homomorphism
$$ \Ind: \bigoplus_{\overset{H\subseteq G}\to{\text{elem.}}}
	  \K(H)
	 @>>> \K(G)
$$
is surjective;
the key point is that
the unit element $1$ of $\K(G)$
lies in the ideal generated by
the $\Ind$-images of $\K(H)$
where $H$ runs over all elementary subgroups of $G$.
Therefore, the proposition follows from
applying Brauer's theorem to the finite group $G/N$.
\qed
\enddemo

\proclaim{Corollary~6}
Let $G$ be a reductive group,
and let $N$ be an open normal subgroup of $G$.
Let $\rho$ be a representation of $G$.
Then there exist a finite list of pairs:
$$ (H_1,\sigma_1),\ldots,(H_s,\sigma_s),
\leqno (*)
$$
where,
for each $i=1,\ldots,s$,
\listitems
\litem{(a)}
$H_i$ is an open subgroup of $G$
with $N \subseteq H_i \subseteq G$,
\litem{(b)}
$\sigma_i$ is an irreducible representation of $H_i$, and in fact,
\litem{(c)}
$\Res^{H_i}_N(\sigma_i)$ is
an irreducible representation of $N$,
\endlistitems\noindent
such that
one has an isomorphism of representations of $G$
of the form
$$ \rho\ \oplus
   \Bigl( \bigoplus_{i=1}^t \Ind^G_{H_i}(\sigma_i) \Bigr)
	\quad\cong\quad
   \Bigl( \bigoplus_{j=t+1}^s \Ind^G_{H_j}(\sigma_j) \Bigr)
\leqno (**)
$$
for some $t$ with $1\le t\le s$.
\endproclaim

\remark{Remark}
If one takes $N$ to be
the identity component $G^0$ of $G$,
then property~(c) asserts that
each $\sigma_i$ is Lie-irreducible.
This is the situation which we shall encounter later in \S4.
\endremark

\demo{Proof of Corollary~6}
Proposition~5 tells us that
we can find a finite list of pairs as in~$(*)$,
such that an isomorphism of the form~$(**)$ holds,
such that properties~(a) and (b) are verified,
and such that each $H_i$ is elementary modulo $N$.
Since each $H_i/N$ is then a nilpotent finite group,
proposition~4 allows us
to replace each $H_i$ by a subgroup containing $N$
and each $\sigma_i$ by an irreducible representation
of the corresponding subgroup,
so that, furthermore, property~(c) is also verified.
This proves the corollary.
\qed
\enddemo

\head{\S3.}
Descent of Representations
\endhead

Let $\Gamma$ be a group,
let $k_0$ be a field of characteristic zero,
and let $k$ be a field extension of $k_0$.
In this section,
we prove two criteria (propositions~7 and~9)
for descending a $k$-representation of $\Gamma$
to a $k_0$-representation.

\proclaim{Proposition~7}
Let $\rho$ be a finite-dimensional $k$-representation of $\Gamma$,
which is absolutely irreducible
(i.e. irreducible over an algebraic closure of $k$).
Assume:
\listitems
\litem{(i)}
$\rho$ is defined over
a finite Galois extension $K$ of $k_0$ in $k$;
\litem{(ii)}
for every $\gamma\in\Gamma$,
the trace $\Tr(\rho(\gamma))$ of $\gamma$
with respect to $\rho$
lies in $k_0$;
\litem{(iii)}
there exists some $\alpha\in k_0$
and some $\gamma_0\in\Gamma$
such that
$\alpha$ is an eigenvalue of multiplicity~one
of $\gamma_0$ with respect to $\rho$.
\endlistitems 
\noindent
Then $\rho$ is defined over $k_0$.
\endproclaim

\demo{Proof of Proposition~7}
By~(i), we may assume that
$\rho$ is given as a $K$-matrix representation of $\Gamma$:
$$ \rho: \Gamma @>>> \GL_r(K),
$$
and we let $\Sigma = \Gal(K/k_0)$ be the finite Galois group.
According to~(iii), we may
choose an eigenvector $v \in {K}^{\oplus r}$
of $\rho(\gamma_0)$ with eigenvalue $\alpha$.
By changing basis,
we may assume that
$v$ is the first basis vectors of ${K}^{\oplus r}$;
thus the matrix $\rho(\gamma_0)$ has the form
$$\pmatrix  \alpha	&  *	&\hdots	&  *	\\
		0	&  *	&\hdots	&  *	\\
	     \vdots	& \vdots&\ddots	& \vdots \\
		0	&  *	&\hdots	&  *	\\
  \endpmatrix.
$$
Each $\sigma\in\Sigma$ defines
a $K$-representation
$$ \sigma\rho : \Gamma @>{\rho}>> \GL_r(K)
				  @>{\GL_r(\sigma)}>> \GL_r(K).
$$
Since $\alpha\in k_0$ is invariant under $\Sigma$,
the matrices $\sigma\rho(\gamma_0)$
also have the same form as $\rho(\gamma_0)$ above;
thus
$v$ is also an eigenvector with eigenvalue $\alpha$
of each $\sigma\rho(\gamma_0)$, $\sigma\in\Sigma$.

Assumption~(ii) and the invariance of $k_0$ under $\Sigma$
gives the equality in $k_0$:
$$ \Tr(\sigma\rho(\gamma)) = \Tr(\rho(\gamma))
   \qquad\text{for any $\sigma\in\Sigma$, any $\gamma\in\Gamma$}.
$$
Therefore, by the trace comparison theorem of Bourbaki
(cf.~\cite{B}~\S12, no.~1, Prop.~3),
the $K$-representations
$\sigma\rho$ of $\Gamma$,
for various $\sigma\in\Sigma$,
are all isomorphic over $K$ to $\rho$.
Choose such isomorphisms over $K$:
$$ a(\sigma) :
	(\sigma\rho, {K}^{\oplus r})
	  \ @>{\cong}>>\  
	(\rho, {K}^{\oplus r}),
		\qquad \sigma \in \Sigma.
$$
Since $\rho$ is absolutely irreducible by hypothesis,
any automorphism of it must be a scalar in $K$.
It follows that
each $a(\sigma) \in \GL_r(K)$ is determined
up to a $K$-scalar multiple.
For any $\sigma,\sigma'\in\Sigma$,
the two different ways of
expressing $\sigma'\sigma\rho$ in terms of $\rho$
then gives
$$ a(\sigma'\sigma)
 = (\text{scalar in $K$})\cdot a(\sigma')\cdot \sigma'a(\sigma).
$$

We shall now rigidify the situation.
For each $\sigma\in\Sigma$,
we have the equality
$$ a(\sigma) \cdot \sigma\rho(\gamma_0)
  = \rho(\gamma_0) \cdot a(\sigma),
$$
and the fact that
$v \in {K}^{\oplus r}$ is an eigenvector
of $\sigma\rho(\gamma_0)$ with eigenvalue $\alpha$;
it follows that
$a(\sigma)\,v \in {K}^{\oplus r}$ is an eigenvector
of $\rho(\gamma_0)$ with eigenvalue $\alpha$.
Thanks to the multiplicity-one hypothesis~(iii) on $\alpha$,
$a(\sigma)\,v$ is necessarily
a $K$-scalar multiple of $v$ itself.
Since we are free to adjust
$a(\sigma) \in \GL_r(K)$
by any $K$-scalar multiple,
we may and do assume that
each $a(\sigma)$ maps $v$ to itself.
Thus the matrices $a(\sigma)$, for $\sigma\in\Sigma$,
have the form
$$\pmatrix      1	&  *	&\hdots	&  *	\\
		0	&  *	&\hdots	&  *	\\
	     \vdots	& \vdots&\ddots	& \vdots \\
		0	&  *	&\hdots	&  *	\\
  \endpmatrix,
$$
and it follows that
the matrices $\sigma'a(\sigma)$, for $\sigma,\sigma'\in\Sigma$,
also have the same form above,
which implies that
each $\sigma'a(\sigma)$ also maps $v$ to itself.
Therefore, we now have
$$ a(\sigma'\sigma) 
 = a(\sigma')\cdot \sigma'a(\sigma)
	\qquad\text{for any $\sigma,\sigma'\in\Sigma$}.
$$
By Hilbert Theorem~90 for $\GL_r$,
there exists some $b\in\GL_r(K)$
such that
$$ a(\sigma) = b\cdot \sigma b^{-1}
	\qquad\text{for each $\sigma\in\Sigma$}.
$$
Using $b^{-1}\in\GL_r(K)$ for a change of basis,
we obtain the $K$-representation
$\tilde\rho := b^{-1}\,\rho\,b$ defined by
$$ \tilde\rho: \Gamma @>>> \GL_r(K),
	\qquad \gamma \mapsto b^{-1}\,\rho(\gamma)\,b,
$$
which is isomorphic over $K$ to $\rho$.
A straightforward computation now shows that
the matrices
$$ \tilde\rho(\gamma) \in\GL_r(K),
	\qquad\text{for $\gamma\in\Gamma$},
$$
are all fixed under the action of the Galois group $\Sigma$;
in other words,
$\sigma\tilde\rho = \tilde\rho$ for any $\sigma\in\Sigma$.
Thus the representation $\tilde\rho$ factorizes as
$$ \Gamma @>>> \GL_r(k_0)
		  \hookrightarrow \GL_r(K).
$$
So $\tilde\rho$ is defined over $k_0$,
and the same is therefore true for $\rho$.
\qed
\enddemo

\proclaim{Lemma~8}
Let $M, N$ be $k_0$-representations of $\Gamma$.
\listitems
\litem{(i)}
The canonical homomorphism of $k$-vector spaces
$$ k\otimes_{k_0} \Hom_{k_0\Gamma}(M,N) 
   @>>>
   \Hom_{k\Gamma}(k\otimes_{k_0}M, k\otimes_{k_0}N ) 
$$
is injective;
it is surjective
if $M$ is finitely generated as a left $k_0\Gamma$-module.
\litem{(ii)}
The canonical homomorphism of $k$-vector spaces
$$ k\otimes_{k_0}\Ext^1_{k_0\Gamma}(M,N)
   @>>>
   \Ext^1_{k\Gamma}(k\otimes_{k_0}M, k\otimes_{k_0}N ) 
$$
is injective
if $M$ is finitely generated as a left $k_0\Gamma$-module.
\endlistitems
\endproclaim

\remark{Remarks}
a) If $M$ is {\it finitely presented\/}
as a left $k_0\Gamma$-module,
the lemma follows from the well-known
``change of rings'' isomorphisms
applied to $k_0\Gamma \hookrightarrow k\Gamma$
(see \cite{R} Th.~2.39 for instance).
Of course, if $M$ is a finite-dimensional
$k_0$-representation of $\Gamma$,
then it is automatically
a finitely generated left $k_0\Gamma$-module;
however,
it need not be finitely presented as a
left $k_0\Gamma$-module.

\noindent
b) When $\Gamma$ is a finite group,
the group ring $k_0\Gamma$ is left-noetherian,
so a finite dimensional $k_0$-representation $M$ of $\Gamma$
is finitely presented as a left $k_0\Gamma$-module,
and the lemma follows from a) above.
But since we will use the lemma
when $\Gamma$ is a profinite group,
and we could not identify a satisfactory reference
for the corresponding result,
we find it prudent to give a complete proof here.

\noindent
c) The proof below actually shows that
the lemma holds in slightly greater generality:
it suffices to assume that
$k_0$ is any commutative ring,
and that $k$ is a $k_0$-algebra
which is {\it free as a $k_0$-module\/}.
\endremark

\demo{Proof of Lemma~8}
We first show that
the canonical homomorphism
$$\matrix\format\r&\ \c\ &\l		\\
  k\otimes_{k_0}\Hom_{k_0\Gamma}(M,N)
  & @>>> &
  \Hom_{k\Gamma}(k\otimes_{k_0}M, k\otimes_{k_0}N )	\\
  \alpha\otimes\phi \quad
  & \mapsto &
  \quad (\beta\otimes m \mapsto \alpha\beta\otimes\phi(m))
\endmatrix
$$
is injective.
Choose a basis $\{ e_i\in k \ :\ i\in I \}$ of $k$
as a $k_0$-vector space.
Then the $k_0\Gamma$-module $k\otimes_{k_0}N$
is the direct sum of
the $k_0\Gamma$-submodules $e_i\otimes N$:
$$ k\otimes_{k_0}N 
   \ \cong\ 
   \bigoplus_{i\in I} e_i\otimes N;
$$
likewise,
the $k_0$-vector space $k\otimes_{k_0}\Hom_{k_0\Gamma}(M,N)$
is the direct sum of
the corresponding $k_0$-subspaces $e_i\otimes\Hom_{k_0\Gamma}(M,N)$:
$$ k\otimes_{k_0}\Hom_{k_0\Gamma}(M,N)
   \ \cong\ 
   \bigoplus_{i\in I} e_i\otimes\Hom_{k_0\Gamma}(M,N).
$$
Any $\phi\in k\otimes_{k_0}\Hom_{k_0\Gamma}(M,N)$
is therefore equal to a sum
$$ \phi \ =\  \sum_{i\in I} e_i\otimes\phi_i
$$
for some uniquely determined
$\phi_i\in\Hom_{k_0\Gamma}(M,N)$, $i\in I$,
all but finitely of which are the zero-map.
Suppose $\phi$ lies in the kernel of the canonical homomorphism.
Then for any $m\in M$,
one has
$$ \sum_{i\in I} e_i\otimes\phi_i(m) \ =\ 0
   \qquad\text{in}\qquad
   k\otimes_{k_0}N 
   \ \cong\ 
   \bigoplus_{i\in I} e_i\otimes N,
$$
so $\phi_i(m) = 0$ in $N$ for each $i\in I$.
It follows that $\phi = 0$,
which is what we want.

If $M$ is finite free as a left $k_0\Gamma$-module,
then it follows from
the functorial properties of $\Hom$ and $\otimes$ that
the canonical homomorphism is an isomorphism.
In general,
if $M$ is finitely generated as a left $k_0\Gamma$-module,
let
$$ 0 @>>> K  @>>> F @>>> M @>>> 0
$$
be a short exact sequence of left $k_0\Gamma$-modules
with $F$ finite free.
Then
$$ 0 @>>> \Hom_{k_0\Gamma}(M,N)
     @>>> \Hom_{k_0\Gamma}(F,N)
     @>>> \Hom_{k_0\Gamma}(K,N)
$$
is an exact sequence of $k_0$-vector spaces.
From this and the fact that $k$ is flat over $k_0$,
we obtain the following commutative diagram with exact columns:
$$\CD
0	@.	0	\\
@VVV		@VVV	\\
k\otimes_{k_0}\Hom_{k_0\Gamma}(M,N)
	@>>>
		\Hom_{k\Gamma}(k\otimes_{k_0}M,k\otimes_{k_0}N)	\\
@VVV		@VVV	\\
k\otimes_{k_0}\Hom_{k_0\Gamma}(F,N)
	@>{\cong}>>
		\Hom_{k\Gamma}(k\otimes_{k_0}F,k\otimes_{k_0}N)	\\
@VVV		@VVV	\\
k\otimes_{k_0}\Hom_{k_0\Gamma}(K,N)
	{\kern10pt\lhook\kern-8pt}\joinrel@>>>
		\Hom_{k\Gamma}(k\otimes_{k_0}K,k\otimes_{k_0}N)	\\
@VVV		@VVV
\endCD
$$
where the middle horizontal arrow is an isomorphism
and the bottom horizontal arrow is injective,
by what we have already shown.
A diagram chase shows that
the top horizontal arrow is surjective.
This proves part~(i).

For part~(ii),
we write down the next terms in the above commutative diagram:
$$\CD
k\otimes_{k_0}\Hom_{k_0\Gamma}(F,N)
	@>{\cong}>>
		\Hom_{k\Gamma}(k\otimes_{k_0}F,k\otimes_{k_0}N)	\\
@VVV		@VVV	\\
k\otimes_{k_0}\Hom_{k_0\Gamma}(K,N)
	{\kern10pt\lhook\kern-8pt}\joinrel@>>>
		\Hom_{k\Gamma}(k\otimes_{k_0}K,k\otimes_{k_0}N)	\\
@VVV		@VVV	\\
k\otimes_{k_0}\Ext^1_{k_0\Gamma}(M,N)
	@>>>
		\Ext^1_{k\Gamma}(k\otimes_{k_0}M,k\otimes_{k_0}N) \\
@VVV		@VVV	\\
0	@.	0
\endCD
$$
By part~(i),
the top horizontal arrow is an isomorphism
and the middle horizontal arrow is injective.
A diagram chase shows that
the bottom horizontal arrow is injective.
This proves part~(ii).
\qed
\enddemo

\proclaim{Proposition~9 (E.~Noether -- M.~Deuring)}
Let $\rho$, $\tau$ and $\pi$ be
semisimple finite-dimensional $k$-representations of $\Gamma$
such that
$$ \rho \oplus \tau \ \cong\  \pi.
$$
Suppose $\tau$ and $\pi$ are defined over $k_0$.
Then $\rho$ is also defined over $k_0$.
\endproclaim

\demo{Proof of Proposition~9}
Our argument here is adapted from that given
for representations of finite groups
(see \cite{H}~Theorem~37.6 for instance).
The proposition is clear when $\tau = 0$.
We proceed by induction on the rank $\rank(\tau)$ of $\tau$;
hence assume that $\rank(\tau) \ge 1$.
By hypothesis,
there exist $k_0$-representations $\tau_0$, $\pi_0$ of $\Gamma$
such that
$$ \tau \cong  k\otimes_{k_0}\tau_0,
	\qquad
   \pi  \cong  k\otimes_{k_0}\pi_0.
$$
For any finite-dimensional
$k_0$-representations $M, N$ of $\Gamma$,
we have the canonical inclusion:
$$ \Ext^1_{k_0\Gamma}(M,N)
   \hookrightarrow
   k\otimes_{k_0}\Ext^1_{k_0\Gamma}(M,N)
   \lhook\joinrel@>>{\text{Lemma~8}}>
   \Ext^1_{k\Gamma}(k\otimes_{k_0}M, k\otimes_{k_0}M);
$$
this fact and the hypothesis that
$\tau$, $\pi$ are semisimple
as $k$-representations of $\Gamma$
imply that
$\tau_0$, $\pi_0$ are semisimple
as $k_0$-representations of $\Gamma$.

Let $\sigma_0 \subseteq \tau_0$ be
an irreducible constituent of
the $k_0$-representation $\tau_0$ of $\Gamma$.
Then
$$ k\otimes_{k_0}\Hom_{k_0\Gamma}(\sigma_0,\pi_0)
   @>{\cong}>{\text{Lemma~8}}>
   \Hom_{k\Gamma}(\sigma,\pi)
   \ \cong\ 
   \Hom_{k\Gamma}(\sigma,\rho \oplus \tau)
$$
contains
$$
   \Hom_{k\Gamma}(\sigma,\tau)
   @<{\cong}<{\text{Lemma~8}}<
   k\otimes_{k_0}\Hom_{k_0\Gamma}(\sigma_0,\tau_0)
   \ \ne 0,
$$
whence
$\Hom_{k_0\Gamma}(\sigma_0,\pi_0) \ne 0$.
Thus $\sigma_0$ is also
an irreducible constituent of
the $k_0$-representation $\pi_0$ of $\Gamma$.
Therefore,
$$ \tau_0 \ \cong\  \tau'_0 \oplus \sigma_0,
	\qquad
    \pi_0 \ \cong\   \pi'_0 \oplus \sigma_0,
$$
for some
semisimple $k_0$-representations $\tau'_0$ and $\pi'_0$
of $\Gamma$.
Letting
$$ \tau' := k\otimes_{k_0}\tau'_0,
   \quad
   \pi'  := k\otimes_{k_0}\pi'_0,
   \quad
   \sigma := k\otimes_{k_0}\sigma_0,
$$
we obtain an isomorphism
$$ \rho \oplus \tau' \oplus \sigma
	\ \cong\ 
                \pi' \oplus \sigma
$$
of semisimple $k$-representations of $\Gamma$,
and hence an equality of their $k$-valued trace functions:
$$ \Tr(\rho(g)) + \Tr(\tau'(g)) + \Tr(\sigma(g))
	\ =\ 
		  \Tr(\pi'(g))  + \Tr(\sigma(g))
  \qquad\text{for every $g\in\Gamma$}.
$$
Applying the trace comparison theorem of Bourbaki
(cf.~\cite{B}~\S12, no.~1, Prop.~3)
to the equality
$$ \Tr(\rho(g)) + \Tr(\tau'(g))
	\ =\ 
		  \Tr(\pi'(g))
  \qquad\text{for every $g\in\Gamma$},
$$
we obtain an isomorphism
$$ \rho \oplus \tau' \ \cong\  \pi'
$$
of semisimple $k$-representations of $\Gamma$.
Since $\rank(\tau') < \rank(\tau)$,
our induction hypothesis shows that
$\rho$ is defined over $k_0$.
\qed
\enddemo

\head{\S4.}
Proof of Main Theorem
\endhead

We shall now prove the main theorem
stated in the introduction.

Thus,
let $\F_q$ be a finite field of characteristic $p$,
let $\ell\ne p$ be a prime number,
let $X$ be a normal variety over $\F_q$,
and let $\Lcal$ be a lisse $\Qellbar$-sheaf on $X$,
which is irreducible,
and whose determinant is of finite order.
Let $E\subset\Qellbar$ denote
the number field given by hypothesis~(1) applied to $(X,\Lcal)$;
thus for every closed point $x$ of $X$,
the polynomial
$$ \det(1-T\,\Frob_x,\Lcal)
$$
has coefficients in $E$.
We may replace the finite field $\F_q$
by its algebraic closure
in the function field $\kappa(X)$ of $X$,
and hence assume that
$X$ is geometrically connected over $\F_q$;
this allows us to use the results in \S1.
Let $\etabar @>>> X$ be a geometric point of $X$,
and set
$$ \Gamma := \pi_1(X,\etabar),
	\qquad
   G := \Garith(\Lcal,\etabar).
$$
Let
$$ \rho_\Lcal : \Gamma @>>> \GL(\Lcal_\etabar)
$$
denote
the monodromy $\Qellbar$-representation of $\Gamma$
corresponding to $\Lcal$,
and let
$$ \rho : G \hookrightarrow \GL(\Lcal_\etabar)
$$
denote
the faithful representation of $\Garith(\Lcal,\etabar)$ on $\Lcal_\etabar$.

By proposition~1~(ii),
$G$ is a
(possibly non-connected) semisimple algebraic group.
We apply corollary~6 to
the representation $\rho$ of $G$,
with $N := G^0 = \Garith(\Lcal,\etabar)^0$,
to obtain
a finite list of pairs as in~$(*)$,
satisfying the properties~(a), (b) and~(c) listed there,
such that an isomorphism of representations of $G$
of the form~$(**)$ holds.

Consider any pair $(H_i,\sigma_i)$ in $(*)$.
By property~(a),
the identity component $H_i^0$ of $H_i$
is a connected semisimple algebraic group
(in fact it is $\Garith(\Lcal,\etabar)^0$),
which is therefore equal to its own commutator subgroup;
hence the $1$-dimensional representation $\det(\sigma_i$) of $H_i$,
given by the determinant of $\sigma_i$,
factors through $H_i/H_i^0$,
and so is given by a character of $H_i$ of finite order.
This and properties~(b) and~(c)
show that
each $\sigma_i$ is
a Lie-irreducible representation of $H_i$,
and its determinant is of finite order.

Set
$$ \Gamma_i := (\rho_\Lcal)^{-1}(H_i) \subseteq \Gamma.
$$
Then $\Gamma_i$ is an open subgroup of $\Gamma$,
corresponding to
a finite etale cover
$ X_i @>>> X $
of $X$ by
a connected variety $X_i$
pointed by the geometric point~$\etabar$;
we identify $\Gamma_i$ with
the arithmetic fundamental group $\pi_1(X_i,\etabar)$ of $X_i$.
If $V_i$ is the representation space of $\sigma_i$,
then the composite homomorphism
$$ \sigma_{\Fcal_i} : 
       \Gamma_i @>{\rho_\Lcal}>> H_i @>{\sigma_i}>> \GL(V_i)
$$
is a $\Qellbar$-representation of $\Gamma_i$
which corresponds to
a lisse $\Qellbar$-sheaf $\Fcal_i$ on the variety $X_i$.
It follows from the corresponding properties of $\sigma_i$ that
$\Fcal_i$ is Lie-irreducible,
and its determinant is of finite order.
By hypothesis~(1) applied to $(X_i, \Fcal_i)$,
there is a number field $E_i \subset \Qellbar$
such that
for every closed point $x$ of $X_i$,
the polynomial
$$ \det(1-T\,\Frob_x, \Fcal_i)
$$
has coefficients in $E_i$;
and by proposition~2,
there is some $\alpha_i\in\Qellbar$
and some closed point $x_0^{(i)}$ of $X_i$
such that
$\alpha_i$ is an eigenvalue of multiplicity~one
of $\Frob_{x_0^{(i)}}$ acting on $\Fcal_i$.
It follows that
$\alpha_i$ is algebraic over the number field $E_i$.

Let 
$$ \rho_{\Lcal_i}
   :=
   \Ind^\Gamma_{\Gamma_i}(\sigma_{\Fcal_i})
$$
be the $\Qellbar$-representation of $\Gamma$
induced from $\sigma_{\Fcal_i}$,
and let
$$ F :=
  \text{composite of $E_1(\alpha_1),\ldots,E_s(\alpha_s)$
	and $E$ in $\Qellbar$}.
$$
It is clear that
$F$ is a finite extension of $E$ in $\Qellbar$.
The isomorphism~$(**)$ implies that
for any closed point $x$ of $X$, one has
$$ \Tr(\rho_{\Lcal}(\Frob_x))  \ +\  
   \sum_{i=1}^t \Tr(\rho_{\Lcal_i}(\Frob_x))
	\ = \ 
   \sum_{j=t+1}^s \Tr(\rho_{\Lcal_j}(\Frob_x))
	\qquad\text{(equality in $F \subset \Qellbar$)}.
\leqno (***)
$$
We shall now show that
the number field $F$
satisfies the conclusion of assertion~(\threeprime).

To that end,
pick a place $\lambda'$ of $F$
lying over a prime number $\ell'\ne p$,
and choose an algebraic closure $\Qellprimebar$ of $F_{\lambda'}$.
By hypothesis~(3) applied to $(X,\Lcal)$ and each $(X_i,\Fcal_i)$,
there exist
irreducible lisse $\Qellprimebar$-sheaves
$\Lcal'$ on $X$ and $\Fcal_i'$ on $X_i$,
which are {\it compatible with $\Lcal$ and $\Fcal_i$\/} respectively;
i.e. for each closed point $x$ of $X$, one has
$$  \det(1-T\,\Frob_x,\Lcal')
  = \det(1-T\,\Frob_x,\Lcal)
	\qquad\text{(equality in $F[T]$)},
\leqno (1)
$$
and for each $i=1,\ldots,s$
and each closed point $x$ of $X_i$, one has
$$  \det(1-T\,\Frob_x,\Fcal_i')
  = \det(1-T\,\Frob_x,\Fcal_i)
	\qquad\text{(equality in $F[T]$)}.
\leqno (2)
$$
It follows that
$\Lcal'$ has the same rank as $\Lcal$
(and each $\Fcal_i'$ has the same rank as $\Fcal_i$).
It also follows from these compatibility relations that
$$ \alpha_i\in F\subset F_{\lambda'}\subset\Qellprimebar
   \quad\text{is an eigenvalue of multiplicity~one
	of $\Frob_{x_0^{(i)}}$ acting on $\Fcal_i'$}.
\leqno (3)
$$

Let $\rho_{\Lcal'}$ denote
the irreducible monodromy $\Qellprimebar$-representation of $\Gamma$
corresponding to $\Lcal'$,
and let $\sigma_{\Fcal_i'}$ denote
the irreducible monodromy $\Qellprimebar$-representation of $\Gamma_i$
corresponding to $\Fcal_i'$.
Let 
$$ \rho_{\Lcal_i'}
   :=
   \Ind^\Gamma_{\Gamma_i}(\sigma_{\Fcal_i'})
$$
be the $\Qellprimebar$-representation of $\Gamma$
induced from $\sigma_{\Fcal_i'}$.
From~(1) and~(2),
we deduce that
for each closed point $x$ of $X$, one has
$$  \Tr(\rho_{\Lcal'}(\Frob_x))
  = \Tr(\rho_{\Lcal}(\Frob_x))
	\qquad\text{(equality in $F$)},
\leqno (4)
$$
and for each $i=1,\ldots,s$
and each closed point $x$ of $X_i$, one has
$$  \Tr(\sigma_{\Fcal_i'}(\Frob_x))
  = \Tr(\sigma_{\Fcal_i}(\Frob_x))
	\qquad\text{(equality in $F$)},
\leqno (5)
$$
whence for each $i=1,\ldots,s$
and each closed point $x$ of $X$, one has
$$  \Tr(\rho_{\Lcal_i'}(\Frob_x))
  = \Tr(\rho_{\Lcal_i}(\Frob_x))
	\qquad\text{(equality in $F$)}.
\leqno (6)
$$
Combining the equalities $(4),(6)$ with~$(***)$,
we see that
for any closed point $x$ of $X$, one has
$$ \Tr(\rho_{\Lcal'}(\Frob_x))  \ +\ 
   \sum_{i=1}^t \Tr(\rho_{\Lcal_i'}(\Frob_x))
	\ = \ 
   \sum_{j=t+1}^s \Tr(\rho_{\Lcal_j'}(\Frob_x))
	\qquad\text{(equality in $F \subset \Qellprimebar$)}.
\leqno (***')
$$
By \u{C}ebotarev's density theorem,
this equality of traces,
as an equality in $\Qellprimebar$,
holds for every element of $\Gamma$.
Therefore, by the trace comparison theorem of Bourbaki
(cf.~\cite{B}~\S12, no.~1, Prop.~3),
we obtain an isomorphism of
semisimple $\Qellprimebar$-representations of $\Gamma$:
$$ \rho_{\Lcal'}\ \oplus
   \Bigl( \bigoplus_{i=1}^t \rho_{\Lcal_i'} \Bigr)
	\quad\cong\quad
   \Bigl( \bigoplus_{j=t+1}^s \rho_{\Lcal_j'} \Bigr).
\leqno (7)
$$

Consider the (absolutely) irreducible $\Qellprimebar$-representation
$\sigma_{\Fcal_i'}$ of $\Gamma_i$.
We wish to apply proposition~7 to this representation;
so let us check that the hypotheses there are verified.
\listitems
\litem{(i)}
By the definition of lisse $\Qellprimebar$-sheaves
(cf. \cite{D}~(1.1.1) ---
 alternatively, apply \cite{KSa}~Remark~9.0.7),
the $\Qellprimebar$-representation $\sigma_{\Fcal_i'}$
is defined over
a finite extension of $\Qellprime$,
which we may of course assume to be
finite Galois over $F_{\lambda'}$.
\litem{(ii)}
From~(5),
we see that
for every closed point $x$ of $X_i$,
the trace $\Tr(\sigma_{\Fcal_i'}(\Frob_x))$
of $\Frob_x\subset\Gamma_i$ with respect to $\sigma_{\Fcal_i'}$
lies in $F_{\lambda'}$;
so it follows from \u{C}ebotarev's density theorem that
the trace $\Tr(\sigma_{\Fcal_i'}(\gamma))$
of every element $\gamma\in\Gamma_i$ with respect to $\sigma_{\Fcal_i'}$
lies in $F_{\lambda'}$.
\litem{(iii)}
Finally, from~(3),
we know that
$\alpha_i\in F_{\lambda'}$
is an eigenvalue of multiplicity~one
of $\Frob_{x_0^{(i)}}\subset\Gamma_i$
with respect to $\sigma_{\Fcal_i'}$.
\endlistitems
\noindent
Hence proposition~7 shows that
$\sigma_{\Fcal_i'}$ is defined over $F_{\lambda'}$.
Then each $\rho_{\Lcal_i'}$,
being induced from $\sigma_{\Fcal_i'}$,
is also defined over $F_{\lambda'}$.
Therefore, in $(7)$,
the two representations in parentheses
are defined over $F_{\lambda'}$.
Proposition~9 now shows that
$\rho_{\Lcal'}$ is also defined over $F_{\lambda'}$,
and hence the lisse $\Qellprimebar$-sheaf $\Lcal'$
is defined over $F_{\lambda'}$;
in other words,
there exists
a lisse $F_{\lambda'}$-sheaf $\Lcal_{\lambda'}$ on $X$
such that
$\Lcal' \cong \Lcal_{\lambda'}\otimes_{F_{\lambda'}} \Qellprimebar$.
The asserted properties of $\Lcal_{\lambda'}$
follow from this isomorphism and~$(1)$.

This completes the proof of our main theorem.

\enddocument